\documentclass[11pt]{article}

\usepackage{amsmath,amsfonts,amscd, amssymb,theorem}

\newcommand{\doubleheaddownarrow}{\big\downarrow\kern-3.325mm\downarrow}
\newcommand\Oh{\mathcal O}
\newcommand\Z{\mathbb Z}
\newcommand\Q{\mathbb Q}

\newcommand{\iso}{\cong}

\newcommand\Pf{\operatorname{Pf}}

\newcommand{\wD}{\widetilde{D}}

\newcommand{\Grass}{\operatorname{Grass}}

\renewcommand\div{\operatorname{div}}

\newcommand\Spin{\operatorname{Spin}}
\newcommand\Ker{\operatorname{Ker}}

\newcommand\Image{\operatorname{Im}}

\newcommand\Proj{\operatorname{Proj}}
\newcommand\pos{\operatorname{pos}}
\newcommand\negat{\operatorname{neg}}
\newcommand\Ext{\operatorname{Ext}}

\newcommand\Hom{\operatorname{Hom}}

\newcommand{\skewentry}{ \kern -1cm -\mathrm{sym} \kern -1cm}

\newcommand{\unpr}{\operatorname {unpr}}
\newcommand{\Sym}{\operatorname {Sym}}
\renewcommand\P{\operatorname {\mathbb P}}

 \newtheorem{theorem}{Theorem}[section]
 \newtheorem{lemma}[theorem]{Lemma}
 
 \newtheorem{cor}[theorem]{Corollary}
 {
 \theorembodyfont{\rmfamily}
 \newtheorem{defn}[theorem]{Definition}
 \newtheorem{exa}[theorem]{Example}

 \newtheorem{rem}[theorem]{Remark}

 }
 \newenvironment{pf}{\paragraph{Proof}}{\par\medskip}
 
 \newcommand{\qed}{\ifhmode\unskip\nobreak\fi\quad\ensuremath\square}
\newcommand{\QED}{\ifhmode\unskip\nobreak\fi\quad\ensuremath{\mathrm{QED}}}

\numberwithin{equation}{section}

\title{ Towards a general theory of unprojection }
    
\author{Stavros Argyrios Papadakis}
\date {July 2006}

\begin{document}
\maketitle

\begin {abstract}
Unprojection is an effort, initiated by Miles Reid, to develop
an algebraic language for the study of birational geometry.
\cite{Ki} contains motivation and examples, and poses the
problem of developing a general theory of unprojection. The
main purpose of the present work is to suggest a general definition
of unprojection, and to show that it indeed generalizes previous 
work done in the topic. In addition, in 
Section~\ref{sec!reid_suzuki_example} we present an unprojection 
analysis of an example of Reid and K. Suzuki, and
Section~\ref{sec!examples} contains more
examples.
\end{abstract}

\section  {Introduction}

Birational geometry is an old and important field  of
algebraic geometry. Since late 1970s there has been spectacular
progress, especially in the establishment of the Mori minimal model
program  for threefolds due to work of S. Mori and many others.
The methods of the Mori minimal model
program are often abstract and 
cohomological, while explicit birational geometry  \cite{EBG}
initiated by Reid and  A. Corti aims to study the objects 
(such as Fano $3$-folds) and the maps between them (such as 
birational contractions) in more detail
on specific situations.  Unprojection plays an
important role in this study
 and has found many applications, for example in the birational 
geometry of Fano $3$-folds  \cite{CPR} and \cite{CM},  in the construction of
weighted complete  intersection K3 surfaces and  Fano
$3$-folds  \cite{Al}, and  in the study of
Mori flips \cite{BrR}. 

\cite {Ki} discusses more examples and applications of unprojection, and
poses the problem of developing a general theory of unprojection.
The cases that have
been studied so far are the unprojection of type Kustin--Miller (or type I)
\cite{KM}, \cite{PR} and \cite{P}, the generic case of type II
unprojection  \cite{P2}, and the generic case of type III unprojection
\cite{P3}, while \cite{R} contains examples of type IV unprojection.
The definitions of unprojection in \cite{P2} and \cite{P3} apply
only to their respective generic cases, while the definitions in \cite{KM}
and \cite{PR} need strong Gorenstein assumptions.

In Section~\ref{sec!gendefnofunpr} we propose a general definition
of unprojection  (Definition~\ref{defn!defnofunprojnew}),
while in Section~\ref{sec!homomequivextensions} we use the 
well--known general machinery of homological algebra to
write down explicitly some of the constructions 
needed.

In Section~\ref{sec!dfnofidealunp} we study in some detail
the important special case of unprojection of an
ideal, while in Section~\ref{sec!compatibility}
we prove that Definition~\ref{defn!defnofunprojnew} indeed 
generalizes those of \cite{PR}, \cite{P2} and \cite{P3}.

Section~\ref{sec!reid_suzuki_example} presents an 
analysis of a construction of  Reid and Suzuki \cite{RS}
which is an unprojection  but is not the unprojection of 
an ideal. Finally, Section~\ref{sec!examples} contains more
examples of  unprojection analysis of rings appearing
in geometry.

A very interesting open question stated in Remark~\ref{rem!openquestion}
is to study under which conditions 
good properties of the unprojection initial data are preserved by
the unprojection ring. Another open question is whether 
unprojection can be used for an inductive treatment
of families of rings  arising in geometry such as
homogeneous coordinate rings of Grassmannians and other
homogeneous spaces, cf. \cite{KM}~Section~2.

ACKNOWLEDGEMENTS: I wish to thank Winfried Bruns, Miles Reid,
and Frank--Olaf Schreyer for  important discussions and
Margarida Mendes Lopes for comments that improved the
presentation of the material.
Part of this work was financially supported by 
the Deutschen Forschungsgemeinschaft Schr 307/4-2.

\section {General definition of unprojection} 
 \label{sec!gendefnofunpr}

Assume $\Oh_X$ is a commutative ring with unit, $M$
is an $\Oh_X$-module, and 
\[  
    \phi \colon \Oh_X \to M
\]
is a homomorphism of $\Oh_X$-modules.  We assume that
there exists an $\Oh_X$-regular element $q \in \Oh_X$ such that
$q M = 0$.  We will define an $\Oh_X$-algebra 
$\unpr_{\Oh_X} \phi$ which we will call the unprojection algebra
of $\phi $.

Since a map $M \to \Oh_X/(q)$  is a 
homomorphism of $\Oh_X$-modules if and only if
is a homomorphism of $\Oh_X/(q)$-modules,  
by Rees lemma (\cite{BH}~Lemma~3.1.16) there are 
canonical isomorphisms
\[
   \Ext^1_{\Oh_X} (M, \Oh_X)
      \iso  \Hom_{\Oh_X} (M, \Oh_X/(q)),
\]
and
\[
   \Ext^1_{\Oh_X} (\Ext^1_{\Oh_X} (M, \Oh_X), \Oh_X)
      \iso \Hom_{\Oh_X} (\Hom_{\Oh_X} (M, \Oh_X/(q)), \Oh_X/(q)).
\]

In Section~\ref{sec!homomequivextensions} we explicitly 
write down, using the well--known general machinery of homological
algebra, how an extension induces a homomorphism
and also how a homomorphism induces an extension.

As a consequence, composing the natural double dual map
\[
   M \to \Hom_{\Oh_X} (\Hom_{\Oh_X} (M, \Oh_X/(q)), \Oh_X/(q))
\]
with $\phi$ and taking the value of the composition at
$1_{\Oh_X} \in \Oh_X$ we obtain an extension
\begin{equation}  \label {eqn!theextension}
    0 \to \Oh_X  \xrightarrow{\,i \, }  \  Q   \to \Ext^1_{\Oh_X} (M, \Oh_X) \to0.
\end {equation}

Again from general principles of homological algebra, if we choose another 
$\Oh_X$-regular element $q'$ with $q'M =0$ we would
get an extension  isomorphic to (\ref {eqn!theextension}).

We define the
$\Oh_X$-algebra $R_1$ with
\[
       R_1 =   \frac {\Sym_{\Oh_X} Q}  {(i(1_{\Oh_X})-1)},
\]
where $\Sym_{\Oh_X} Q$ is the symmetric algebra of the 
$\Oh_X$-module $Q$. We also  define
 the multiplicatively closed subset  $T \subset \Oh_X $ 
\[
    T = \{ t \in \Oh_X \colon t \text { is both } \Oh_X \text { and }
     M  \text{-regular } \} \subset \Oh_X.
\]

Therefore, the set $T$ consists of the elements of $\Oh_X$  which multiplied
by any nonzero element of $\Oh_X$ or $M$   give a nonzero
element of the respective $\Oh_X$-module. In particular $T$
contains the invertible elements of $\Oh_X$.

\begin {defn}  \label{defn!defnofunprojnew}
The unprojection algebra $\unpr_{\Oh_X} \phi$ 
is the $\Oh_X$-algebra
\[
  \unpr_{\Oh_X}\phi  = R_1/J, 
\]  
where 
\[
    J = \{ u \in R_1 \colon \text { there exists } t \in T \text { with }
    tu = 0 \in R_1 \}.
\] 
\end {defn}

\begin {rem}
An important case for the applications is when a morphism
$\widetilde {\phi} \colon \wD  \to X$
of affine schemes having codimension one image in $X$ is given. 
In this case we set
\[ 
   \phi \colon \Oh_X  \to M =\widetilde{\phi}_* \Oh_{\wD}
\]
to be the homomorphism of $\Oh_X$-modules induced
by the morphism $\widetilde{\phi}$, compare the Reid--Suzuki 
example in Section~\ref{sec!reid_suzuki_example}.
A particular case, which we study 
in some detail in Section~\ref{sec!dfnofidealunp},
is when  $D \subset X$ is a codimension one 
subscheme and $\widetilde {\phi} \colon D \to X$
is the inclusion morphism.
\end {rem}

\section {The relation between extensions and homomorphisms}
  \label{sec!homomequivextensions}

In the following $\Oh_X$ is a commutative ring, $A,B$
are two $\Oh_X$-modules, and $q \in \Oh_X$ is an 
$\Oh_X$ and $B$-regular element such that $qA = 0$. 
By Rees lemma (\cite{BH}~Lemma~3.1.16)
\[ 
   \Ext^1_{\Oh_X} (A,B)
      \iso \Hom_{\Oh_X} (A, B/(q)).
\] 
We use the well--known general machinery of homological algebra to
write down explicitly the correspondence betweeen 
homomorphisms $A \to B/(q)$ and extensions  $0 \to B \to Q \to A \to 0$ 
used in Section~\ref{sec!gendefnofunpr}.

\subsection {Construction of the homomorphism given an extension}
 \label {subs!homomfromextension}

Assume we are given an extension
\begin{equation}  \label {eqn!secondextension}
    0  \to B \xrightarrow{\, q_1 \, }   Q  \xrightarrow{\, q_2 \,}  A \to 0.
\end {equation}
We will define a homomorphism  of $\Oh_X$-modules
\[
   g \colon A \to  B/ (q) 
\]
as follows.

Let $a \in A$. Choose lifting $\widetilde{a} \in Q$. Then
$q\widetilde{a} \in \ker q_2$ (since $qA =0$), so there exists unique
$b \in B$ such that  $q\widetilde{a} =b$ (equality in $Q$).
We set
\[
      g(a)  =  b + (q) \in B/(q).
\] 

Assume $\widetilde{a}' \in Q$ is another lifting of $a$, and 
$b' \in B$ such that  $q\widetilde{a}' =b'$  Since
$  \widetilde{a} - \widetilde{a}' \in \Ker q_2$ and 
(\ref {eqn!secondextension}) is exact,  there exists
$b_3 \in B$ with
\[
       \widetilde{a} - \widetilde{a}' = b_3  
\]
(equality in $Q$). As a consequence
\[
    b - b' = q(\widetilde{a} - \widetilde{a}') = qb_3,
\]
hence
\[
      b + (q) = b' + (q),
\]
(equality in $B/(q)$)
and therefore $g$ is well defined, independent of the choice of the
lifting of $a$.

\subsection {Construction of the extension given a homomorphism} 
    \label {subs!extensionfromhomom}

Assume now we are given a homomorphism 
\[
     g \colon A \to B/(q) 
\]
We will use $g$ to define an extension 
\[
    0  \to B \xrightarrow{\, q_1 \, }   Q  \xrightarrow{\, q_2 \,}  A \to0.
\]
of $\Oh_X$-modules.

Fix a generating set $a_i, i \in I$ for $A$, denote by 
$F$ the free $\Oh_X$-module with basis $e_i, i\in I$, 
and by 
\[ 
   p_1 \colon F \to  A
\]
the surjective $\Oh_X$-module homomorphism 
with $p_1 (e_i) = a_i$.

Moreover, fix (arbitrary) set--theoretic lifting
\[
   \widetilde{g}  \colon A \to B 
\]
of $g$.

\begin {lemma}  \label {lemm!liftingcompatibility}
Assume $r_i \in \Oh_X$ with $i \in I$ 
and all except a finite number of $r_i$ equal to $0$, 
such that $\sum_{i \in I} r_i a_i = 0$ in $A$.
We then have
\[
   \sum_{i \in I} r_i \widetilde{g} (a_i) \in (q) \subset B.
\]
(The map $\widetilde{g}$ in general is not a homomorphism
of $\Oh_X$-modules, so it may happen that
 $\; \sum_{i \in I} r_i \widetilde{g} (a_i) \not= 0 \in B$.)
\end {lemma}

\begin {pf}
   Clear, since $\sum_{i \in I} r_i g(a_i) = 0 \in B/(q)$
and $\widetilde{g}$ is a lifting of $g$.
\QED \medskip
\end {pf}

By Lemma~\ref{lemm!liftingcompatibility}, if  
$\sum_{i \in I} r_i a_i = 0$ in $A$, there exists unique 
(since $q$ is $\Oh_X$-regular) $b \in B$  such that
\[
    qb = \sum_{i \in I} r_i \widetilde{g} (a_i)
\]
in $B$.  
We define an $\Oh_X$-submodule $M \subset B \times F$ with
\[
    M = \{ (b, \sum_{i \in I} r_i e_i) \in B \times F\colon
       \sum_{i \in I} r_i a_i = 0 \in A \; \text { and } \;
          qb = \sum_{i \in I} r_i \widetilde{g} (a_i) \in B
       \}
\]
and we set 
\[
     Q =  (B \times F) / M.
\]

\begin {lemma} 
We have an extension of  $\Oh_X$-modules,
\[
    0  \to B \xrightarrow{\, q_1 \, }   Q  \xrightarrow{\, q_2 \,}  A \to0.
\]
where $q_1$ is the map
\[
  b \mapsto  (b,0) + M \in Q
\]
   and  $q_2$ is the map 
\[
  (b, \sum_{i \in I} r_i e_i) + M \mapsto \sum_{i \in I} r_i a_i \in A.
\]
\end{lemma}

\begin {pf}
The map $q_2$ is well defined, since  if $(b, \sum_{i \in I} r_i e_i) \in M$
we have $\sum_{i \in I} r_i a_i = 0 \in A$.

We also have that $q_1$ is injective. Indeed, assume $(b,0) \in M$. This
implies that $qb = 0 \in B$, and since $q$ is $B$-regular we have $b=0$.

It is clear that 
\[
   q_2 \circ  q_1 = 0.
\]
Assume now  $(b, \sum_{i \in I} r_i e_i) + M  \in \ker q_2$.
This implies that $\sum_{i \in I} r_i a_i = 0 \in A$. By
Lemma~\ref{lemm!liftingcompatibility} we have
\[
   \sum_{i \in I} r_i \widetilde{g} (a_i) \in (q) \subset B, 
\] 
so there exists  $b' \in B$ such that
\[
   qb' = qb  -  \sum_{i \in I} r_i \widetilde{g} (a_i)
\]
(equality in $B$). As a result
\[
   (b-b', \sum_{i \in I} r_i e_i)  \in M,
\]
hence 
\[
  (b, \sum_{i \in I} r_i e_i) + M =(b',0) + M  = q_1(b') \in \Image q_1.
\]

It is clear that $q_2$ is surjective, which finishes the proof of 
the lemma.
\QED \medskip
\end {pf}

It is easy to check that the extension does not depend  
on the choices of the
generating set $a_i, i \in I$ of $A$ and of the lifting $\widetilde{g}$
of $g$.

\section {Unprojection of an ideal} \label{sec!dfnofidealunp}

In the following, $\Oh_X$ is a commutative 
ring with unit, and $I_D \subset \Oh_X$ an ideal 
containing an $\Oh_X$-regular element $q \in I_D$.
We set $\Oh_D = \Oh_X/I_D$  and denote by
\[
   \phi \colon \Oh_X \to \Oh_D
\]
the natural projection map.

\begin {defn}  \label{defn!defnofidealunproj}
The unprojection algebra $\unpr_{\Oh_X}I_D$ 
of the ideal $I_D \subset \Oh_X$ 
is the $\Oh_X$-algebra
\[
    \unpr_{\Oh_X}I_D = \unpr_{\Oh_X} \phi,
\]  
where $\unpr_{\Oh_X} \phi$ has been defined in 
Definition~\ref{defn!defnofunprojnew}.
\end{defn}

It is easy to see that the extension
corresponding to the projection map $\phi$ is just
the natural short exact sequence 
\begin {equation} \label{eqn!adjunctionseq}
  0 \to \Oh_X \to \Hom_{\Oh_X}(I_D, \Oh_X) \to 
     \Ext^1_{\Oh_X} (\Oh_D, \Oh_X) \to 0
\end {equation}
obtained by 
applying the derived functor of $\Hom_{\Oh_X}(-, \Oh_X)$ to the
exact sequence $0\to I_D\to\Oh_X\to\Oh_D\to0$
(cf. \cite{PR} p.~563).

In particular, $Q = \Hom_{\Oh_X}(I_D, \Oh_X)$, and
hence we have  that
\[
    \unpr_{\Oh_X}I_D = R_1/ J,
\]  
where 
\[
   R = \Sym_{\Oh_X} \Hom_{\Oh_X}(I_D, \Oh_X),
\] 
$i \in \Hom_{\Oh_X}(I_D, \Oh_X)$ is the inclusion  map
$I_D \to \Oh_X$, $R_1 =  R / (i-1)$, 
\[
    J = \{ u \in R_1 \colon \text { there exists } t \in T \text { with }
    tu = 0 \in R_1 \},
\]
and $T \subset \Oh_X $ is the multiplicatively closed subset 
\[
    T = \{ t \in \Oh_X \colon t \text { is both } \Oh_X \text { and }
          \Oh_D \text{-regular } \} \subset \Oh_X.
\]

\begin {rem}
In general, the ring $\unpr_{\Oh_X}I_D$ is not graded - for an
important exception see Remark~\ref{rem!grading} below. However,   
the natural grading of $R$ induces an increasing filtration
\[
    F_0 \subset F_1 \subset F_2 \subset \ldots 
\] 
of  $\unpr_{\Oh_X}I_D$, where $F_k$ is the image of the
direct sum of the first $k+1$ graded components of 
$R$ under the natural
projection map $R \to \unpr_{\Oh_X}I_D$.
\end{rem}

\begin{rem} \label{rem!grading}
Let  $\Oh_X$ be $\Z$-graded and $I_D \subset \Oh_X$
a homogeneous ideal. 
We call  an $\Oh_X$-module homomorphism 
$\widetilde{s} \colon I_D \to \Oh_X$ graded, if there exists
$k \in \Z$,  called the degree of $\widetilde{s}$, such that $\widetilde{s}$ sends
homogeneous elements of $I_D$ of degree $a$ to homogeneous
elements of $\Oh_X$ of degree $a+k$ for all $a \in \Z$.

Assume that $\Hom_{\Oh_X}(I_D, \Oh_X)$ is generated by graded 
homomorphisms, a sufficient condition for that is 
that $\Oh_X$ is Noetherian (cf. \cite{BH} p.~33). Then,  there is a 
unique natural grading on $R$ extending the gradings
of $\Oh_X$ and $\Hom_{\Oh_X}(I_D, \Oh_X)$. Under this
grading, $i$ is homogeneous of degree $0$, hence the ideal
$(i-1) \subset R $ is homogeneous and $R_1$ and $\unpr_{\Oh_X}I_D$
become graded rings. 

\end{rem}

\subsection {Relation with $\Oh_X[I_D^{-1}]$}   \label{subssec!sectionrelns}

From now on, we assume that $I_D$ 
contains an $\Oh_X$-regular element.  We denote
by $K(X)$ the total quotient ring of $\Oh_X$, i.e.,
the localization of $\Oh_X$ with respect to the multiplicatively
closed subset of $\Oh_X$ consisting of all $\Oh_X$-regular elements,
and we set
\[
   I_D^{-1} = \{ f \in K(X) \colon f I_D \subset \Oh_X \}.
\]

$I_D^{-1}$ is an $\Oh_X$-submodule of $K(X)$, which under the 
above assumption that $I_D$ contains an $\Oh_X$-regular element
is  naturally isomorphic to $\Hom_{\Oh_X}(I_D, \Oh_X)$.
Indeed, the map 
\[
  I_D^{-1} \to  \Hom_{\Oh_X}(I_D, \Oh_X)
\]  
sends 
$f \in I_D^{-1}$ to the multiplication map by $f$, while the inverse
map sends $g \in \Hom_{\Oh_X}(I_D, \Oh_X)$ to $g(q)/q \in I_D^{-1}$,
where $q \in I_D$ is any $\Oh_X$-regular element. In the following
we will identify these two isomorphic $\Oh_X$-modules.

We denote by  $\Oh_X [I_D^{-1}]$  the  $\Oh_X$-subalgebra of $K(X)$ 
generated by $I_D^{-1}$. 
The inclusion $I_D^{-1} \subset \Oh_X [I_D^{-1}]$
and the universal property of $R=\Sym_{\Oh_X}  I_D^{-1}$ induce a
surjective homomorphism of $\Oh_X$-algebras
\[
       \rho \colon R \to \Oh_X [I_D^{-1}],
\]
which restricted to the first graded part of $R$ is the identity.

\begin {lemma}
The natural inclusion $\Oh_X \subset R$ as degree $0$ graded part
induces injective maps $\Oh_X \to R_1$ and $\Oh_X \to \unpr_{\Oh_X}I_D$.
\end{lemma}

\begin {pf} 
The second claim follows from the first, since $T$ contains only
$\Oh_X$-regular elements. We will show that $\Oh_X \cap(i-1) = 0 \subset R$,
this will prove the first claim.

Let $a \in \Oh_X \cap(i-1)$. There exist homogeneous elements $b_0, \dots ,b_r$
of $R$, with degree of $b_t$ equal to $t$, such that
\[
    a = (1-i)(b_0 + \dots +b_r).
\]
Comparing homogeneous degrees we get $0=b_ri = ai^{r+1}$. Using the map $\rho$
we get $0 = \rho(ai^{r+1})=a$.
\QED \medskip
\end{pf}

Clearly $(i-1) \subset \ker \rho$, so there is an induced map
\[
       \rho^1 \colon R_1 \to \Oh_X [I_D^{-1}].
\]

\begin{lemma}
 We have 
\[
         \ker \rho^1 = J_2,
\]

where 
\[
   J_2 = \{ u \in R_1 \colon \text {there exists } 
   t \in \Oh_X  \text { which is } \Oh_X 
     \text{-regular with } tu = 0 \in R_1 \}.   
\]
\end{lemma}

\begin {pf}

 It is clear that 
 $J_2 \subset \ker \rho^1$. We will show
 the opposite inclusion. By the assumptions, there
 exists $q \in I_D$ which is $\Oh_X$-regular. Therefore,
 if $a \in I_D^{-1}$, there exists 
 $z \in \Oh_X \subset R_1$ (depending on $a$)  with
 \[
   q a - z = 0   \in R_1.
 \]

 Let $u \in \ker \rho^1$. Since  
 $R_1$ is generated as a ring by  $I_D^{-1}$,
 $u$ is a polynomial in elements of $I_D^{-1}$, hence 
 for a sufficiently large integer $n$ there
 exists $z \in \Oh_X \subset R_1$ with  
 \[
   q^n u - z = 0   \in R_1.
 \]
Since  $u \in \ker \rho^1$ we necessarily have $\rho^1 (z) = 0$,
hence $z = 0$.
\QED \medskip
\end{pf}

By the above, the map $\rho^1$ factors through the natural quotient 
maps
\[
   R_1 \to   \unpr_{\Oh_X}I_D \to  \Oh_X[I_D^{-1}]
\]
The last map $\unpr_{\Oh_X}I_D \to  \Oh_X[I_D^{-1}]$ is
often (examples of Section~\ref{sec!examples}, normal case of 
unprojection of type Kustin--Miller 
(\cite {PR} Remark~1.3), generic type II \cite{P2}, 
generic type III \cite{P3}) but not always an isomorphism. 

\begin {exa} 
Consider the Kustin--Miller unprojection pair  
(cf. \cite{P} Section~4)
\[
     I_D = (x,y) \subset \Oh_X = k[x,y]/(x^2-y^3)
\]
and  set
\[
   u = \frac{y^2}{x} = \frac{x}{y} \in I_D^{-1}.
\] 
Then, $u^2-y$ is zero in $\Oh_X[I_D^{-1}] \subset K(X)$, 
but, using Theorem~\ref{thm!generaloftypeI} below, 
$u^2-y$ is nonzero in $\unpr_{\Oh_X}I_D$. 
\end{exa}

\begin {rem}  \label{rem!onlylinears}
 It is obvious that when $\rho^1$ is an
isomorphism, we have that $\unpr_{\Oh_X}I_D$ is
isomorphic to $\Oh_X[I_D^{-1}]$ and to $R_1$ (case of only
linear relations). In the generic type II unprojection \cite{P2}, 
$\rho^1$ is not an isomorphism (there exist
quadratic relations) but still $\unpr_{\Oh_X}I_D$ is
isomorphic to $\Oh_X[I_D^{-1}]$.
\end{rem}

\begin {rem} In our experience, we have often found easier to 
first study  the ring $\Oh_X[I_D^{-1}]$ and then use the results
to study the unprojection ring $\unpr_{\Oh_X}I_D$.
\end{rem}

\begin {rem} \label{rem!adjunctsequence}
Unlike the case where $\Oh_X$ is Gorenstein, the $\Oh_X$-module
$\Ext^1_{\Oh_X}(\Oh_D,\Oh_X)$  which appears in 
(\ref{eqn!adjunctionseq})
in general depends both on $\Oh_D$ and
$\Oh_X$. For an example, fix
\[
    \Oh_D = \frac{k[x_1,x_2,x_3]} {(x_1,x_2,x_3)},
\]
where $k$ is any field, and set 
\[
    (\Oh_X)_n = \frac{k[x_1,x_2,x_3]} {(I_X)_n},
\]
where $(I_X)_n$ is generated by the maximal minors of a 
general $n \times n+1 $ 
matrix with linear entries in $x_i$.
\end {rem}

\begin {rem} \label{rem!openquestion}
 It will be very interesting to study under which conditions 
good properties of the unprojection initial data are preserved by
the unprojection ring. Compare \cite{PR}~Theorem~1.5, 
\cite{P2}~Theorem~2.15 and \cite{P3}~Theorem~3.5.
\end{rem}

\section {Compatibility with previously defined unprojections} \label{sec!compatibility}

\subsection {Case of unprojection of type Kustin--Miller}

In this subsection, we prove that Definition~\ref{defn!defnofidealunproj}
generalizes Definition~1.2 of \cite{PR}.

Assume we are under the assumptions of \cite{PR}~Section~1.
That is, $\Oh_X$ is a local Gorenstein ring, and $I_D=(f_1, \dots ,f_r)
 \subset \Oh_X$
is a codimension one ideal such that the quotient ring $\Oh_D=\Oh_X/I_D$
is also Gorenstein. We choose  as in \cite{PR} an injective map 
$\widetilde{s} \colon I_D \to \Oh_X$
generating the $\Oh_X$-module $\Hom_{\Oh_X}(I_D, \Oh_X)/\Oh_X$, and we set
$g_i = \widetilde{s}(f_i) \in \Oh_X$.

Let
\[
    \psi \colon \Oh_X[S] \to  R
\]
be the $\Oh_X$-algebra homomorphism with $\psi(S) =\widetilde{s} \in R_1$, where  $\Oh_X[S]$
is the polynomial ring over $\Oh_X$ in one variable. The map $\psi$ is
not surjective, but induces two surjective maps
\[
   \psi^1 \colon \Oh_X[S] \to  R_1 
\]   
and
\[
   \psi^2 \colon \Oh_X[S] \to   \unpr_{\Oh_X}I_D. 
\]

The following lemma gives a presentation of
$I_D^{-1} = \Hom_{\Oh_X}(I_D, \Oh_X)$ as $\Oh_X$-module.

\begin {lemma} As an $\Oh_X$-module 
\[
    I_D^{-1}  \iso  \frac {\Oh_X s_0 \oplus \Oh_X s_1} {(f_is_1 - g_is_0)},
\]    
where $s_0$ corresponds to the inclusion $i\colon I_D \to \Oh_X$, and $s_1$
corresponds to $\widetilde{s}$.
\end{lemma}

\begin{pf}  Assume $l_1,l_2 \in \Oh_X$ with 
\[
      l_1 i  - l_2 \widetilde{s} = 0 \in I_D^{-1}.
\]
Then 
\[
   l_2 \widetilde{s} = 0 \in  \Hom_{\Oh_X}(I_D, \Oh_X)/\Oh_X \iso \Oh_D,
\]   
hence $l_2 \in I_D$. Fix $q \in I_D$ an $\Oh_X$-regular element, such
an element  exists since $\Oh_X$ is Cohen--Macaulay and $I_D$ 
has codimension one. Then
\[
   l_1 q = (l_1i) (q) = l_2 \widetilde{s}(q) =\widetilde{s}(l_2q)
 = q \widetilde{s}(l_2),
\]
therefore  $l_1 = \widetilde{s}(l_2)$.
\QED \medskip
\end{pf}

Using elementary properties of the symmetric algebra of a module (cf.
\cite{Ei} p.~570~Prop.~A2.2 )
we have the following corollary.

\begin {cor}  \label{cor!kerofq1}
We have that
\[
    \ker \psi^1 = (Sf_i - g_i).
\]
\end{cor}

\begin {theorem}   \label{thm!generaloftypeI}
We have that
\[
    \ker \psi^2 = (Sf_i - g_i).
\]
As a consequence, Definition~\ref{defn!defnofidealunproj} of unprojection 
generalizes Definition~1.2 of \cite{PR} and 
\[        
    \unpr_{\Oh_X}I_D  \iso  R_1.
\]
\end{theorem}

\begin{pf}
  Let  $h=h(S) = a_nS^n + \dots + a_0 \in \ker \psi^2$, with $a_i \in \Oh_X$ and
  $a_n \not= 0$.  We prove by induction on the degree $n$ of $h$ that 
  $h \in   (Sf_i - g_i)$.  
  
  Assume $n=0$. Then $h \in \Oh_X$ and $th = 0 \in \Oh_X$ for $t \in T$ implies
  that $h=0$, since $T$ contains only $\Oh_X$-regular elements.
 
 Assume $n \geq 1$, and that the result is true for all polynomials of
 degree strictly less than $n$. Using  Corollary~\ref{cor!kerofq1}, and the
 definition of $\unpr_{\Oh_X}I_D$, there exist $l_i(S) \in \Oh_X[S]$ and 
 $t \in T$, with
 \[
    th(S) = \sum l_i(S) (Sf_i-g_i)  \in \Oh_X[S].
 \]
 Write
 \[
     l_i = m_{k_i} S^{k_i} + \text { lower terms },
 \]
 with $ m_{k_i} \in \Oh_X $ nonzero. 

We remark that for $u_i \in \Oh_X$, 
 $\sum u_if_i =0$ implies $\widetilde{s}(\sum u_if_i) =0$, so $\sum u_ig_i = 0$.
 Therefore, $ta_n \in (f_i)=I_D$.
 By definition, $T$ contains only $\Oh_D$-regular elements, hence $a_n \in (f_i)$.
 As a consequence, there exists $h' \in \Oh_X[S]$ with degree strictly less
 than $n$, with $h-h' \in (Sf_i-g_i)$. By the inductive hypothesis 
 $h' \in (Sf_i-g_i)$, so also $h \in (Sf_i-g_i)$ which proves the theorem.
 \QED \medskip
\end{pf}

\subsection {Case of generic type II unprojection}

In this subsection, we prove that Definition~\ref{defn!defnofidealunproj}
generalizes Definition~2.2 of \cite{P2}. 

We use the notations of \cite{P2} Section~2. In addition, we define
$J_L \subset \Oh_X[T_0, \dots ,T_k]$ to be the ideal generated by 
all affine linear polynomials $f_{i,j,p}^a$ and $f_{j,p}^b$.

\begin {lemma}  \label{lem!compatwithtypeII}
Using the notations of \cite{P2} Section~2, we have 
that  
\[
     a_{11}g_{ij}^a   \in J_L
\]
for all $i,j$ with $i+j \leq k$, and that
\[
     a_{11}g_{ij}^b   \in J_L
\]
for all $i,j$ with $i+j \geq k+1$.
\end{lemma}

\begin {pf}  We prove the first inclusion, the second follows
by similar arguments. For $r=1,2, \dots ,i$ we define 
elements $y_r, u_r \in \Oh_X[T_0,\dots ,T_k]$ by
\[
     y_r = (T_j (a_{r,1}T_{i-r+1}+a_{r+1,1}T_{i-r}) -
          T_0 (a_{r,1}T_{i+j-r+1} + a_{r+1,1}T_{i+j-r}))(-1)^{r+1} 
\]
and
\[
    u_r = (T_j f_{r,1,i-r}^a - T_0f_{r,1,i+j-r}^a)(-1)^{r+1}.
\]
From the identity
\[
    a_{11} (T_iT_j - T_0T_{i+j}) = y_1 + y_2 + \dots + y_i,
\]
it follows that the element 
\[
    a_{11}g_{ij}^a  -u_1 - \dots -u_i
\]
is affine linear in $T_p$ and in the kernel of $\phi$. As a 
consequence, it is an element of $J_L$.
 \QED \medskip
\end{pf}

From the above  Lemma~\ref{lem!compatwithtypeII} and 
Proposition~2.6 of \cite{P2}, it follows that 
Definition~\ref{defn!defnofidealunproj}
generalizes Definition~2.2 of \cite{P2}.

\subsection {Case of generic type III unprojection}

In this subsection, we prove that Definition~\ref{defn!defnofidealunproj}
generalizes Definition~3.3 of \cite{P3}.

We use the notations of \cite{P3} Section~3.

From \cite{P3}~Theorem~3.5, it follows  that the natural map
$R \to \Oh_X[I_D^{-1}]$ induces an isomorphism $R/(i-1) \iso \Oh_X[I_D^{-1}]$.
Therefore, using Remark~\ref{rem!onlylinears}, we have that
\[
    R_1 \iso \unpr_{\Oh_X}I_D \iso \Oh_X[I_D^{-1}]  .   
\]

\section {The Reid--Suzuki example}  \label{sec!reid_suzuki_example}

The aim of this section is to  analyse using unprojection a construction
from \cite{RS}~p.~235, which provides an example of unprojection which is 
not the unprojection of an ideal.

We fix a field $K$, the projective line $\P^1$ with 
homogeneous coordinates $v,w$ over $K$,
and three points  $P_1 = [1,0] = (w=0), P_2 = [0,1] = (v=0)$ and 
$Q=[1,1]= (v-w =0)$ of $\P^1$.

\begin {rem} 
Recall that for  a $\Q$-divisor $D$ of $\P^1$,   $H^0(\P^1,D)$ is
the $K$-vector space
\[
   H^0(\P^1,D) = \{ f \in K(\P^1)^* \colon \div (f) + D \text { is effective }
  \Q\text {-divisor } \} \cup \{ 0 \}.
\]
\end{rem}
Consider now the two $\Q$-divisors  
\[
    D_1 = \frac{3}{5}P_1 + \frac{4}{5}P_2 -Q, \quad \; \;
        D_2 = \frac{3}{5}P_1 + \frac{4}{5}P_2 
\]
of $\P^1$,    and the corrresponding graded rings
\begin {eqnarray*}
      \Oh_X & = & \bigoplus_{k \geq 0} H^0(\P^1, kD_1), \\
      \Oh_Y & = &  \bigoplus_{k \geq 0} H^0(\P^1, kD_2).
\end {eqnarray*}

 By \cite{Wa}
$\Oh_X$ and $\Oh_Y$ are both Cohen--Macaulay but not Gorenstein.
  
\subsection {Explicit calculations for the rings $\Oh_X$ and $\Oh_Y$}
  \label {subs!explicitcalcul}

Define $y,z,t,u_1,u_2 \in K(\P^1)$ with 
\begin{eqnarray*}
     y &=& \frac{(v-w)^2}{vw},\; \; z = \frac{(v-w)^3}{v^2w}, \;
     \; t = \frac{(v-w)^4}{v^3w},  \\
     u_1 &=& \frac{(v-w)^5}{v^4w},\; \; u_2 = \frac{(v-w)^5}{v^2w^3}.
\end{eqnarray*}

We consider these elements as homogeneous
elements of $\Oh_X$ of degrees $2,3,4,5,5$ respectively.
It is easy to see that they generate the
$K$-algebra $\Oh_X$, and
\begin{equation} 
  \Oh_X  \iso  K[Y,Z,T,U_1,U_2] / I,
\end{equation} 
with
\[
   I = (Z^2-YT,h_1= ZT-YU_1,h_2=Y^4-ZU_2,T^2-ZU_1,Y^3Z-TU_2,Y^3T-U_1U_2)
\] 
where $Y,Z,T,U_1, U_2$ are indeterminants of degrees $2,3,4,5,5$ 
respectively.

\begin {rem}  \label{rem!m2codeforR1}
 Macaulay 2 code to produce  $\Oh_X$
\begin {verbatim}
kk = ZZ/101
S = kk[v,w]
R = kk[y,z,t,u1,u2, Degrees => {2,3,4,5,5}]
M = matrix {{ v*w, v*w^2,v*w^3, v*w^4, v^3*w^2}}
I = kernel (map (S, R, M))
--  we have \Oh_X = S/I
-- codim I --answer 3 -- res I --answer  r, r^6, r^8 , r^3,  
-- \Oh_X  is CM, but not Gorenstein
\end{verbatim}
\end{rem}

For the ring $\Oh_Y$, we take coordinates $s_0,s_1,y,z,t,u_1,u_2$
where $y,z,t,u_1,u_2$ are as above and $s_0, s_1 \in K(\P^1)$ with
\begin{equation}
    s_0 = 1, \; \; s_1 = \frac {v-w}{w}.
\end{equation}
We consider $s_0$ and $s_1$ as homogeneous elements of $\Oh_Y$ of degrees
1 and 2 respectively, and we also notice that 
\[
   s_0  = \frac{zy}{y^2-t}, \; \;  s_1 = \frac{y^3}{y^2-t}.
\] 

Using a computer algebra system such as Singular \cite{GPS01}
or Macaulay 2 \cite{GS93-08},  we easily get 
\begin{equation} 
  \Oh_Y  \iso  K[S_0,S_1,Y,Z,T,U_1,U_2] / J,
\end{equation} 
where $S_0,S_1,Y,Z,T,U_1,U_2$ are indeterminants of degrees 
$1,2,2,3,4,5,5$ respectively,
and
\[
   J = J_1 + J_2 + J_3 
\] 
with $J_1 = I$, 
\begin{eqnarray*}
J_2 = &  & (-ZS_0+YS_1-Y^2, -TS_0+ZS_1-YZ,  (Y^2-T)S_0-YZ, \\
      & &   h_3=  -U_1S_0+TS_1-YT, (YZ-U_1)S_0-YT, \\
      & &  h_4=(U_1-U_2)S_0+Y^3+YT,  h_5= YTS_0-U_1S_1, \\
      & &    Z(U_1-U_2)S_0+ (Y^2U_1+TU_2), T(U_1-U_2)S_0+  (YZU_1+ U_1U_2), \\
      & & h_6= -Y(Y^2+T)S_0 + (U_1+U_2)S_1 -YU_2, \\
      & & h_7 = 2YU_2S_0-Y(Y^2+T)S_1-ZU_2),
\end{eqnarray*}
(the linear relations between $s_0$ and $s_1$) and
\[
    J_3  =  (h_8= f_4=YS_0^2-S_1^2 +ZS_0+ Y^2)
\]
(the quadratic relation).

\begin {rem}  \label{rem!m2codeforR2}
Macaulay 2 code to produce  $\Oh_Y$ 
\begin {verbatim}
kk = ZZ/101
S = kk[v,w]
newR = kk[y,z,t,u1,u2,s0,s1, Degrees =>{2,3,4,5,5,1,2}]
newM = matrix {{ v*w*(v-w)^2, v*w^2*(v-w)^3, v*w^3*(v-w)^4, 
      v*w^4*(v-w)^5, v^3*w^2*(v-w)^5, v*w, (v-w)*v^2*w}}
newI = kernel (phi2= map (S, newR, newM))
--  we have \Oh_Y = newR / newI
-- codim newI  --answer 5 
-- res newI --answer  r, r^15, r^40, r^45, r^24, r^5   
-- \Oh_Y  is CM  but not Gorenstein
\end{verbatim}
\end{rem}

Define the prime ideal $I_D\subset \Oh_X$, with 
\begin {eqnarray*}
   I_D & = & (z^2-y^3, t-y^2, u_1-yz, u_2-yz)
\end {eqnarray*}
It is easy to see that $I_D$ is the ideal of the 
point $[1,1,1,1,1] \in X = \Proj \Oh_X \subset  \P(2,3,4,5,5) $.

We have 
\[
   \Oh_Y \subset \Oh_X[I_D^{-1}],
\]
however the inclusion is strict.
Indeed, using the following Macaulay 2 code 
continuing the code in Remark \ref{rem!m2codeforR1}

\begin {verbatim}
   ID = ideal(t-y^2, z^2-y^3,u1-y*z, u2-y*z)
   betti presentation Hom(ID, R^1/I)
   presentation Hom(ID, R^1/I)
\end{verbatim}     
we get that $\Hom_{\Oh_X}(I_D, \Oh_X) $ is generated by
\[
   q_0 = \frac{z}{y^2-t}, \; \;  	q_1 = \frac{y^2}{y^2-t},
    \;  \;     q_2 = 1
\]
the first of homogenous degree $-1$, the other two of
homogeneous degrees $0$. Clearly
\[
   s_0 = yq_0 \quad \text { and }  \quad  s_1 =  y q_1.
\]

Another calculation shows that
\begin{eqnarray*}
 (u_1+u_2)h_8 =  & &(s_1-y)h_1+(-2s_0)h_2 +(-z)h_3+(s_0y)h_4 \\
           & & +(-2y)h_5+(-s_1-y)h_6+(s_0)h_7   
\end{eqnarray*}
and it is easy to see that $u_1+u_2$ is both $\Oh_X$ 
and $\Oh_D = \Oh_X/I_D$-regular element. As a consequence,
the quadratic relation $h_8$ multiplied by a regular element
is inside the ideal generated by the linear relations (cf.
Definition~\ref{defn!defnofunprojnew}).

\subsection {The $\Oh_X$-algebra $\Oh_Y$ as unprojection}

Define $\wD = \Proj K[a] \iso \P(1) $ where $a$ is an indeterminate
of degree $1$ and the morphism of schemes
\[
   \widetilde{\phi} \colon  \wD  \to X \subset \P(2,3,4,5,5)
\]
with
\[
    \widetilde{\phi} ([a]) =[a^2, a^3, a^4, a^5, a^5]
\]
induced  by the graded homomorphism of graded $K$-algebras
\[
   \phi \colon  \Oh_X \to K[a]
\]
specified by
\[ 
   \phi(y) =a^2, \; \phi(z)= a^3, \; \phi(t)= a^4, \; 
   \phi(u_1)= \phi(u_2) = a^5.
\]

It is clear that the scheme--theoretic image of $\widetilde{\phi}$ 
is $D \subset X$,
and also that $K[a]$ needs two generators, corresponding to $1$ and $a$,
when viewed as $\Oh_X$-module via $\phi$.

The interpretation of the calculations in 
Subsection~\ref{subs!explicitcalcul} is that we have
\[
       \Oh_Y \iso \unpr_{\Oh_X} \phi,
\]
while $\Oh_Y$ is not isomorphic to $\unpr_{\Oh_X} I_D$.

\section {More examples}  \label{sec!examples}

We discuss below more examples of unprojections related to geometry. 
In all cases we start from a triple of embedded schemes
\[
    D \subset X \subset \P
\]
and construct by unprojection a new embedded scheme 
\[
    Y \subset \P'.
\]
$\P$ and $\P'$ are usual projective spaces, in each case clear from
the construction.   By $I_X, I_Y$ and $I_D$ we will denote the homogeneous
ideals of the corresponding schemes, and by $\Oh_X, \Oh_Y$ etc. 
the corresponding homogeneous  coordinate rings. By abuse of notation,
we will also denote by $I_D$ the homogeneous ideal of $\Oh_X$ corresponding
to $D \subset X$.   In all cases 
\[
    \Oh_Y  = \unpr_{\Oh_X}I_D
\]
and moreover the map $\rho^1$ defined in Section~\ref{subssec!sectionrelns}
is an isomorphism, hence by Remark~\ref{rem!onlylinears}
\[
    \Oh_Y \iso R_1 \iso \Oh_X [I_D^{-1}].
\]

If the minimal resolution of $\Oh_X$ as $\Oh_{\P}$-module taking no accounts
of the twists is

\[
  \Oh_{\P}^{h_0} \xleftarrow{} \Oh_{\P}^{h_1} \xleftarrow{} \dots  \xleftarrow{}
     \Oh_{\P}^{h_d} \xleftarrow{} 0 
\] 
we say that the Betti vector of $X$ is $(h_0,h_1, \dots ,h_d)$, and
similarly for $\Oh_D$ and $\Oh_Y$.

We have checked the calculations of all examples using the computer
algebra programs Macaulay 2 \cite{GS93-08}  and Singular \cite {GPS01}
over the field $\Q$ and the finite field $\Z_{101}$. We believe that 
with a little extra effort one should be able to prove the results 
with coefficients over any field, and even more generally.

\subsection {Rolling factors format}   \label{subsec!rollingfactors}

According to \cite{Ki} Example~10.8, the rolling factors
format first appeared - in a disguised codimension four 
form - in D. Dicks Ph.D. thesis \cite{D}. It has also appeared
in \cite{R1}, \cite{S} and \cite{BCP}.

Let $M$ be the $2 \times 2$ matrix
\[
    M=  \begin {pmatrix}  
         x_1 & x_2 & x_3 & x_4 \\
         y_1 & y_2 & y_3 & y_4
       \end{pmatrix} 
\] 
and denote by $I_M$ the ideal generated by the 
$2 \times 2$ minors of $M$. Let 
\begin {eqnarray*}
   f_1 = a_1x_1 + \dots + a_4x_4  \\
   f_2 = b_1x_1 + \dots + b_4x_4  \\
   g_1 = a_1y_1 + \dots + a_4y_4 \\
   g_2 = b_1y_1 + \dots + b_4y_4 
\end{eqnarray*}
and set 
\[
   I_Y = I_M + (f_1,f_2, g_1, g_2).
\]
The ideal $I_Y$ defines a codimension five
projectively Gorenstein subscheme $Y \subset \P^{15}$
with Betti vector $(1,10,19,19,10,1)$.

Let $N$ be the submatrix of $M$ obtained by
deleting the first column of $M$, and denote by $I_N$
the ideal generated by the $2 \times 2$ minors of $N$. 

Define $X \subset \P^{13}$ (with coordinates of $\P^{13}$ those 
of $\P^{15}$ minus $x_1,y_1$)  by 
\[
   I_X = I_N + ( b_1 f_1 - a_1 f_2,  b_1 g_1 - a_1 g_2).
\]
$X$ is a codimension three projectively Gorenstein
subscheme with Betti vector $(1,5,5,1)$.

Consider  the subscheme $D \subset X$ with 
\[
   I_D = I_X + (a_1,b_1).
\]
The ideal $I_D$ is codimension four projectively 
Cohen--Macaulay with Betti vector $(1,5,9,7,2)$ (actually $I_D$ is a
hyperplane section of $I_N$), and we have 
\[
    \Oh_Y = \unpr_{\Oh_X}I_D  \iso R_1 \iso \Oh_X [I_D^{-1}].
\]

\subsection {Veronese surface}      \label{subsec!veronesesurface}

Let $Y \subset \P^5$ be the Veronese surface. It is well known that if
\[
     M =    \begin {pmatrix}  
         a_{11} & a_{12}  & a_{13} \\
                & a_{22}  & a_{23}  \\
           \text{sym}     &   & a_{33}
       \end{pmatrix} 
\]
is the generic $3 \times 3$ symmetric matrix on the coordinates of $\P^5$, 
then  $I_Y$ is generated by the six $2 \times 2$  minors of $M$. In addition,
$Y$ has codimension three and Betti vector $(1,6,8,3)$. Therefore,
it is projectively Cohen--Macaulay, but not projectively Gorenstein.

Let $N$ be the submatrix of $M$ obtained by
deleting the first row of $M$.

Define $X \subset \P^4$ (with coordinates of $\P^4$ those of $\P^5$ minus
$a_{11}$) with $I_X$ generated by the three $2 \times 2$ minors of $N$.
$X$ has  codimension two and Betti vector $(1,3,2)$. Consider 
the codimension three complete intersection $D \subset \P^4$, with
\[
   I_D = (a_{22},a_{23},a_{33}).
\]

We have   $D \subset X$ and 
\[
    \Oh_Y = \unpr_{\Oh_X}I_D  \iso R_1 \iso \Oh_X [I_D^{-1}].
\]

\subsection {$\Grass(2,6)$}  \label{subsec!Grass}

Let $Y \subset \P^{14}$ be the Grassmannian $\Grass(2,6)$ of dimension $2$ 
vector subspaces of a $6$ dimensional vector space in its natural 
Pl\"ucker embedding. It is well known that if $M = [a_{ij}]$ is the generic 
$6 \times 6$ skew--symmetric matrix on the coordinates of $P^{14}$, then
$I_Y$ is the ideal generated by the $15$ $4 \times 4$ Pfaffians of
$M$ (that is, the Pfaffians of all $4 \times 4$ skew--symmetric 
submatrices of $M$ obtained by deleting from $M$ two rows and the
corresponding two columns). In addition, $Y$ is codimension six and 
projectively Gorenstein, with Betti vector  $(1,15,35,42,35,15,1)$.

Let $N$ be the $5 \times 5$ submatrix of $M$ obtained by  deleting
the first row and the first column of $M$, and let $\Pf(N)$ be the ideal
generated by the $5$ $4 \times 4$ Pfaffians of $N$. 

Define $X \subset \P^{12}$ (with coordinates of $\P^{12}$ those of $\P^{14}$ minus
$a_{12},a_{13}$)  by 
\[
    I_X =  \Pf(N) + ( a_{14} a_{56}- a_{15} a_{46}+a_{16} a_{45} ),
\]
and $D \subset X$ with 
\[
    I_D = I_X + (a_{45}, a_{46}, a_{56}).
\]   

$X$ is codimension four projectively Gorenstein with Betti vector
$(1,6,10,6,1)$, and $D$ is codimension five projectively 
Cohen--Macaulay with Betti vector $(1,6,14,16,9,2)$. We have  
\[
    \Oh_Y = \unpr_{\Oh_X}I_D  \iso R_1 \iso \Oh_X [I_D^{-1}].
\]

\subsection {Spinor Variety}  \label{subsec!spinor}

Denote by $Y \subset P^{15}$ the ten-dimensional spinor 
variety  $\Spin (5,10)$ in its canonical embedding.
It has Betti vector $(1,10,16,16,10,1)$, hence it is 
projectively Gorenstein.

According to $\cite{Ki}$ Problem 8.7,
 \[
    I_Y = (\pos_1, \dots , \pos_5, \negat_1 ,\dots ,\negat_5),
 \]
 where 
 \begin {eqnarray*}
 \pos_1 &  = &  z  y_1 -  x_{23}x_{45}+   x_{24}x_{35}-   x_{25}x_{34} \\
 \pos_2 &  = &  z  y_2 -  x_{13}x_{45}+   x_{14}x_{35}-   x_{15}x_{34} \\
 \pos_3 & = &  z  y_3 -  x_{12}x_{45}+   x_{14}x_{25}-   x_{15}x_{24} \\
 \pos_4 & = &  z  y_4 -  x_{12}x_{35}+   x_{13}x_{25}-   x_{15}x_{23} \\
 \pos_5 & = &  z  y_5 -  x_{12}x_{34} +  x_{13}x_{24} -  x_{14}x_{23} \\
 \negat_1 & = &  x_{12}y_2 -  x_{13}y_3 +   x_{14}y_4 -   x_{15}y_5  \\
 \negat_2 & = &  x_{12}y_1 -  x_{23}y_3 +   x_{24}y_4 -   x_{25}y_5  \\
 \negat_3 & = &  x_{13}y_1 -  x_{23}y_2 +   x_{34}y_4 -   x_{35}y_5  \\
 \negat_4 & = &  x_{14}y_1 -  x_{24}y_2 +   x_{34}y_3 -   x_{45}y_5  \\ 
 \negat_5 & = &  x_{15}y_1 -  x_{25}y_2 +   x_{35}y_3 -   x_{45}y_4.
\end {eqnarray*}

Define $X \subset \P^{14}$ (with coordinates of $\P^{14}$ those of 
$\P^{15}$ minus $z$), by
\[
   I_X = (\negat_1, \dots ,\negat_5).
\]
$X$ is a codimension four almost complete intersection (hence not 
projectively Gorenstein) with  Betti vector $(1,5, 12, 10, 2)$, as 
a consequence it is projectively Cohen--Macaulay.

Consider  the codimension five complete intersection  $D \subset \P^{14}$ 
with 
\[
   I_D = (y_1, y_2, \dots ,y_5).
\]

We have $D \subset X$, and 
\[
    \Oh_Y = \unpr_{\Oh_X}I_D  \iso R_1 \iso \Oh_X [I_D^{-1}].
\]

\bigskip
\noindent
Departamento de Matem\'{a}tica/CAMGSD \\
Instituto Superior T\'{e}cnico  \\
Av. Rovisco Pais  \\ 
1049-001 Lisboa   \\
Portugal    \\
e-mail: spapad@maths.warwick.ac.uk

\end{document}